\documentclass[12pt,a4paper,leqno]{article}

\usepackage{latexsym}
\usepackage{amssymb,exscale}
\usepackage[centertags]{amsmath}
\usepackage{amsthm}
\usepackage[all]{xy}

\usepackage[dvips]{hyperref}

\numberwithin{equation}{section}

\swapnumbers
\theoremstyle{definition}
\newtheorem{theorem}[equation]{Theorem}
\newtheorem*{thm}{Theorem}
\newtheorem{lemma}[equation]{Lemma}
\newtheorem{proposition}[equation]{Proposition}
\newtheorem{corollary}[equation]{Corollary}
\newtheorem{definition}[equation]{Definition}

\newtheorem{counterexample}[equation]{Counterexample}

\renewcommand{\phi}{\varphi}

\newcommand{\D}{\mathrm{d}}

\newcommand{\ti}{\tilde}

\renewcommand{\(}{\bigl(}
\renewcommand{\)}{\bigr)\vphantom{)}}

\newcommand{\imp}{$ \Longrightarrow $ }

\newcommand{\DCS}{\operatorname{DCS}}

\newcommand{\mes}{\operatorname{mes}}

\newcommand{\modO}{{\operatorname{mod}\,0}}
\newcommand{\modI}{{\operatorname{mod}\,1}}

\newcommand{\One}{\mathbf1}

\newcommand{\eps}{\varepsilon}
\newcommand{\si}{\sigma}
\newcommand{\ga}{\gamma}
\newcommand{\om}{\omega}
\newcommand{\Om}{\Omega}
\newcommand{\de}{\delta}

\newcommand{\al}{\alpha}
\newcommand{\be}{\beta}
\newcommand{\Ec}{\mathcal E}
\newcommand{\F}{\mathcal F}

\newcommand{\A}{\mathcal A}

\renewcommand{\Pr}[1]{\mathbb{P}\mskip1.5mu\(\mskip1.5mu#1\mskip1.5mu\)}
\newcommand{\R}{\mathbb R}

\newcommand{\cP}[2]{\mathbb{P}\mskip1.5mu\(\mskip1.5mu#1\mskip1.5mu
 \big|\mskip1.5mu#2\mskip1.5mu\)}

\newcommand{\sif}{$\sigma$\nobreakdash-field}

\newcommand{\normal}[1]{$#1$\nobreakdash-\hspace{0pt}normal}
\newcommand{\set}[1]{$#1$\nobreakdash-\hspace{0pt}set}
\newcommand{\almost}[1]{$#1$\nobreakdash-\hspace{0pt}almost}

\newcommand{\measurability}[1]{$#1$\nobreakdash-\hspace{0pt}measurability}
\newcommand{\measurable}[1]{$#1$\nobreakdash-\hspace{0pt}measurable}

\begin{document}

\title{Random dense countable sets:\\characterization by independence} 

\author{Boris Tsirelson}

\date{}
\maketitle

\begin{abstract}
A random dense countable set is characterized (in distribution) by
independence and stationarity. Two examples are \emph{Brownian local
minima} and \emph{unordered infinite sample}. They are identically
distributed; the former \emph{ad hoc} proof of this fact is now superseded by
a general result.
\end{abstract}

\section*{Introduction}
Random dense countable sets arise naturally from various probabilistic
models. Their examination is impeded by the singular nature of the set $
\DCS(0,1) $ of all dense countable subsets of (say) the interval $ (0,1)
$. This set is not a  Polish space, not even a standard Borel
space. Nevertheless the idea of random elements of $ \DCS(0,1) $ and their
distributions can be formalized. An appropriate framework proposed in
\cite[Sect.~1]{I} is used here.

Two examples of random dense countable sets are compared in \cite{I}. One
example, `Brownian local minima', is the random set
\[
M = \{ s \in (0,1) : \exists \eps>0 \;\> \forall t \in (s-\eps,s) \cup
(s,s+\eps) \;\> B_s < B_t \}
\]
of local minimizers on $ (0,1) $ of the Brownian motion $ (B_t)_t $. The other
example, `unordered infinite sample', is the random set
\[
S = \{ U_1, U_2, \dots \} = \{ s \in (0,1) : \exists n \;\> U_n = s \}
\]
where $ U_1, U_2, \dots $ are independent random variables distributed
uniformly on $ (0,1) $. The main result of \cite{I} states that $ M $ and $ S
$ are identically distributed, which means existence of such a joining between
the Brownian motion $ (B_t)_t $ and the sequence $ (U_n)_n $ that $ M = S $
a.s.

Independence of Brownian increments on disjoint time intervals $ (a,b) $ and $
(c,d) $ implies independence of `fragments' $ M \cap (a,b) $ and $ M \cap
(c,d) $ of $ M $ (see \ref{1.1} for the definition). Independence of $ S \cap
(a,b) $ and $ S \cap (c,d) $ is less evident but true \cite[2.2]{I}. The same
holds for any number of fragments.

Stationarity of Brownian increments should imply stationarity of $ M
$. However, time shifts do not preserve the time interval $ (0,1) $. We have
two options: either replace $ (0,1) $ with the whole $ \R $, or replace linear
shifts $ t \mapsto t+s $ of $ \R $ with cyclic shifts $ t \mapsto t+s \; \modI
$ of $ (0,1) $; I choose the latter option. The cyclic shift is nothing but an
interval exchange transformation: $ (0,1-s) \to (s,1) $ and $ (1-s,1) \to
(0,s) $. Brownian increments on $ (0,1-s) $ and $ (s,1) $ are distributed
identically; taking independence into account we get cyclic stationarity of $
M $ (see \ref{1.9} for the definition).

Cyclic stationarity of $ S $ is evident.

Thus, the main result of \cite{I} is a special case of the following. (See
Definitions \ref{1.2}, \ref{1.9} and Theorem \ref{1.11}.)

\begin{thm}
All random dense countable subsets of $ (0,1) $, satisfying the conditions of
independence and cyclic stationarity, are identically distributed.
\end{thm}

Waiving stationarity we get some other distributions; see Counterexample
\ref{1.6} and Theorem \ref{1.8}.

\section[]{\raggedright Definitions and claims}
\label{sect1}\emph{Throughout, either by assumption or by construction, all probability
spaces are standard.} Recall that a standard probability space (known also as
a Lebesgue-Rokhlin space) is a probability space isomorphic $ (\modO) $ to an
interval with the Lebesgue measure, a finite or countable collection of atoms,
or a combination of both.

According to \cite[Sect.~1]{I}, the set $ \DCS(0,1) $ is of the form $ B/E $,
the quotient set of a standard Borel space $ B $ by an equivalence relation $
E \subset B \times B $. Namely, $ B = (0,1)^\infty_{\ne} $ consists of all
sequences $ (u_n)_n $ of pairwise different points of $ (0,1) $, and $ E $
consists of all pairs $ \( (u_n)_n, (u_{\si(n)})_n \) $ where $ (u_n)_n $ runs
over $ B $ and $ \si $ runs over all permutations of $ \{1,2,\dots\} $.

A map $ X : \Om \to B/E $ is called measurable \cite[Def.~1.3]{I} if it admits
a measurable lifting $ Y : \Om \to B $:
\[
\xymatrix{
 \Om \ar[r]^{Y} \ar[dr]_{X} & B \ar[d]^{\text{canonical projection}}
\\
 & B/E
}
\]

By a random dense countable subset of $ (0,1) $ we mean a measurable map $ \Om
\to \DCS(0,1) $ (or rather an equivalence class of such maps). Its measurable
lifting $ Y = (Y_1,Y_2,\dots) $ may be called also a measurable
\emph{enumeration} of $ X $; $ X = \{ Y_1, Y_2, \dots \} $ in the sense that $
X(\om) = \{ Y_1(\om), Y_2(\om), \dots \} $ for almost all $ \om $.

 Of course,
another interval may be used instead of $ (0,1) $.

Two measurable maps $ X_1 : \Om_1 \to B/E $, $ X_2 : \Om_2 \to B/E $ are
called identically distributed \cite[Def.~1.4]{I} if they can be matched by a
joining $ J $ between $ \Om_1 $ and $ \Om_2 $ (in other words, a measure with
given marginals on $ \Om_1 \times \Om_2 $):
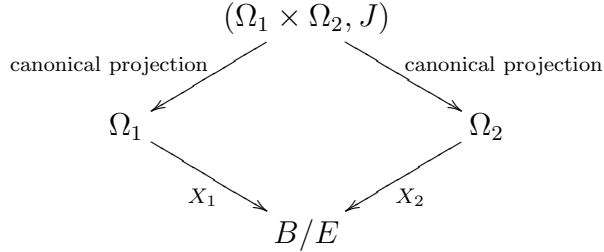
\[
\xymatrix{
 & (\Om_1 \times \Om_2, J)
  \ar[dl]_(.6){\text{canonical projection}\quad}
  \ar[dr]^(.6){\quad\text{canonical projection}}
\\
 \Om_1 \ar[dr]_{X_1} & & \Om_2 \ar[dl]^{X_2}
\\
 & B/E
}
\]

\begin{definition}\label{1.1}
Two measurable maps $ X_1 : \Om \to B/E $, $ X_2 : \Om \to B/E $ are
\emph{independent,} if they admit independent measurable liftings $ Y_1, Y_2 :
\Om \to B $.
\end{definition}

Independence of more than two maps is defined similarly.

\begin{definition}\label{1.2}
Let $ X $ be a random dense countable subset of $ (0,1) $ such that
\begin{equation}\label{1.3}
\Pr{ t \in X } = 0 \quad \text{for each } t \in (0,1) \, .
\end{equation}
We say that $ X $ satisfies the
\emph{independence condition,} if for all $ n = 1,2,\dots $ and all $ 0 = t_0
< t_1 < \dots < t_n = 1 $ the random dense countable subsets $ X_k $ of $
(t_{k-1},t_k) $ defined by $ X_k = X \cap (t_{k-1},t_k) $ are independent.
\end{definition}

As was noted in Introduction, the independence condition is satisfied both for
$ M $ (Brownian local minima) and $ S $ (unordered infinite sample).

\begin{lemma}\label{1.4}
For every random dense countable subset of $ (0,1) $ there exists a measure $
\mu $ on $ (0,1) $ such that
\begin{equation}\label{1.5}
X \cap A = \emptyset \;\; \text{a.s.} \quad \text{if and only if} \quad \mu(A)
= 0
\end{equation}
for all Borel sets $ A \subset (0,1) $.
\end{lemma}

\begin{proof}
A measure lifting $ Y : \Om \to (0,1)^\infty_{\ne} $ enumerates the random
dense countable set $ X $ by random variables $ Y_1, Y_2, \dots $ The measure
$ \mu $ defined by
\[
\mu(A) = \sum_n \frac1{n^2} \Pr{ Y_n \in A }
\]
fits evidently.
\end{proof}

The measure $ \mu $ is determined by $ X $ up to equivalence (mutual absolute
continuity). For $ M $ (Brownian local minima) and $ S $ (unordered infinite
sample) we may take $ \mu = \mes $ (the Lebesgue measure). In general $ \mu $
is nonatomic, otherwise arbitrary.

It may seem that all random dense countable subsets of $ (0,1) $ satisfying
the independence condition and having $ \mu = \mes $ are identically
distributed. However, they are not.

\begin{counterexample}\label{1.6}
We choose a (nonrandom) dense open set $ G \subset (0,1) $ such that $ \mes G
< 1 $; its complement $ C = (0,1) \setminus G $ is a nowhere dense compact set
of positive measure. We take the usual Poisson random subset $ P $ of $ (0,1)
$ (whose intensity measure is the Lebesgue measure) and combine it with $ S $
(unordered infinite sample, independent of $ P $):
\[
X = ( P \cap C ) \cup ( S \cap G ) \, .
\]
It is easy to see that $ X $ is a random dense countable subset of $ (0,1) $
satisfying \eqref{1.3}, the independence condition, and \eqref{1.5} for $ \mu =
\mes $. However, $ X $ and $ S $ are not identically distributed.
\end{counterexample}

In order to exclude such cases we introduce a condition (stronger than
\eqref{1.5} for $ \mu = \mes $):
\begin{equation}\label{1.7}
\begin{aligned}
& X \cap A = \emptyset \text{ a.s.} & & \text{if } \mes(A) = 0 \, , \\
& X \cap A \text{ is infinite a.s.} & & \text{if } \mes(A) > 0
\end{aligned}
\end{equation}
for all Borel sets $ A \subset (0,1) $.

\begin{theorem}\label{1.8}
All random dense countable subsets of $ (0,1) $ satisfying \eqref{1.7} and the
independence condition are identically distributed.
\end{theorem}

The proof is given in Sect.~\ref{sect5}.

We turn to stationarity, that is, invariance under cyclic shifts $ T_s : (0,1)
\to (0,1) $ defined for $ s \in (0,1) $ by
\[
T_s(t) = t+s \; \modI = \begin{cases}
 t+s &\text{for $ t \in (0,1-s) $}, \\
 t+s-1 &\text{for $ t \in (1-s,1) $};
\end{cases}
\]
$ T_s $ is undefined at $ 1-s $, which does not matter as long as our random
sets satisfy \eqref{1.3}.

Given a random dense countable set $ X $ satisfying \eqref{1.3} and a number $
s \in (0,1) $, we get another random dense countable set $ T_s(X) = \{ T_s(t)
: t \in X \} $. Moreover, we may randomize $ s $ as follows. We multiply the
given probability space by the interval $ (0,1) $ (equipped with Lebesgue
measure) and define on the new probability space a measurable function $ Y :
\Om \times (0,1) \to \DCS(0,1) $ by $ Y(\om,s) = T_s (X(\om)) $. Clearly, $ Y
$ is also a random dense countable set.

\begin{definition}\label{1.9}
A random dense countable subset $ X $ of $ (0,1) $ is \emph{stationary} if it
satisfies \eqref{1.3} and $ X $, $ Y $ are identically distributed where $ Y $
is constructed from $ X $ by the random shift, as described above.
\end{definition}

Two evident examples are $ M $ (Brownian local minima) and $ S $ (unordered
infinite sample).

\begin{proposition}\label{1.10}
If a random dense countable subset of $ (0,1) $ is stationary then it
satisfies \eqref{1.7}.
\end{proposition}

The proof is given in Sect.~\ref{sect5}.

\begin{theorem}\label{1.11}
All stationary random dense countable subsets of $ (0,1) $ satisfying the
independence condition are identically distributed.
\end{theorem}

\begin{proof}
Follows immediately from \ref{1.10} and \ref{1.8}.
\end{proof}

\section[]{\raggedright Existence of measures with given marginals}
\label{sect2}Random sets do not appear at all in this section.

Denote by $ M $ the set of all positive Borel measures $ m $ on the square $
(0,1) \times (0,1) $ such that
\[
m \( U \times (0,1) \) \le \mes(U) \, , \quad
m \( (0,1) \times V \) \le \mes(V)
\]
for all Borel sets $ U,V \subset (0,1) $. In other words, both marginals of $
m $ are bounded by the Lebesgue measure.

For every Borel subset $ W $ of the square we define two numbers,
\begin{gather*}
\al(W) = \sup \{ m(W) : m \in M \} \, , \\
\be(W) = \inf \{ \mes(U) + \mes(V) : \( U \times (0,1) \) \cup \( (0,1) \times
 V \) \supset W \} \, ;
\end{gather*}
the infimum is taken over Borel sets $ U,V \subset (0,1) $ satisfying the
indicated condition. In other words, $ W $ does not intersect the product of
their complements, $ \( (0,1) \setminus U \) \times \( (0,1) \setminus V \)
$.

A well-known Strassen's theorem states that $ \al(W) = \be(W) $ for all closed
sets $ W $. The same holds for all \set{F_\si}s, as is shown below.

\begin{lemma}\label{2.1}
If $ W_n \uparrow W $ (that is, $ W_1 \subset W_2 \subset \dots $ and $ W_1
\cup W_2 \cup \dots = W $; of course, these are Borel subsets of the squuare)
then $ \al(W_n) \uparrow \al(W) $.
\end{lemma}

\begin{proof}
Follows immediately from the fact that $ m(W_n) \uparrow m(W) $ for $ m \in M
$.
\end{proof}

\begin{proposition}\label{2.2}
If $ W_n \uparrow W $ then $ \be(W_n) \uparrow \be(W) $.
\end{proposition}

The proof is given after Prop.~\ref{2.4}.

It follows immediately that $ \al(W) = \be(W) $ for all \set{F_\si}s $ W $.

Following \cite[Sect.~5]{I} we may avoid topological notions (closed sets,
\set{F_\si}s). To this end we denote by $ \A $ the algebra of subsets of $
(0,1) \times (0,1) $ generated by all product sets $ U \times V $ where $ U,V
$ are Borel subsets of $ (0,1) $. We consider the class $ \A_\de $ of all sets
of the form $ W_1 \cap W_2 \cap \dots $ where $ W_1, W_2, \dots \in \A $. (All
closed sets belong to $ \A_\de $.) Further, we consider the class $
\A_{\de\si} $ of all sets of the form $ W_1 \cup W_2 \cup \dots $ where $ W_1,
W_2, \dots \in \A_\de $. (All \set{F_\si}s belong to $ \A_{\de\si} $.)

The equality $ \al(W) = \be(W) $ holds for all $ W \in \A_\de $
\cite[5.6]{I}. Therefore (by \ref{2.1}, \ref{2.2}) it holds for all $ W \in
\A_{\de\si} $.

The square $ (0,1) \times (0,1) $ can be replaced with the product $ B_1
\times B_2 $ of two standard Borel spaces (these are Borel isomorphic to $
(0,1) $) equipped with probability measures. Moreover, we may start with the
product $ (\Om,\F,P) $ of two probability spaces $ (\Om_1,\F_1,P_1) $, $
(\Om_2,\F_2,P_2) $ and consider the \sif\ $ \F $ modulo (the
$\si$\nobreakdash-ideal of all) sets of the form $ (A_1 \times \Om_2) \cup
(\Om_1 \times A_2) $ where $ P_1(A_1) = 0 $, $ P_2(A_2) = 0 $ (and their
subsets). All results of this section can be reformulated readily to this more
general framework, together with their proofs.

\begin{lemma}\label{2.3}
Every sequence of pairs of measurable subsets of $ (0,1) $ contains a
subsequence $ (A_n,B_n)_n $ such that the limits
\begin{gather*}
f(x) = \lim_{n\to\infty} \frac{ \#\{k\le n: x\in A_k\} }{ n } \, , \quad
 g(y) = \lim_{n\to\infty} \frac{ \#\{k\le n: y\in B_k\} }{ n } \, , \\
h(x,y) = \lim_{n\to\infty} \frac{ \#\{k\le n: x\in A_k, y\in B_k\} }{ n }
\end{gather*}
exist for almost all $ x,y \in (0,1) $ and
\[
h(x,y) = f(x) g(y)
\]
for almost all $ x,y \in (0,1) $.
\end{lemma}

\begin{proof}
A result of Aldous \cite{A77} (applied to the products) gives us a
subsequence $ (A_n,B_n)_n $ and a measurable map $ P =
(P_{00},P_{01},P_{10},P_{11}) $ from the square $ (0,1) \times (0,1) $ to the
set of probability measures on the four-point set $ \{0,1\} \times \{0,1\} $
such that for almost every point $ (x,y) \in (0,1) \times (0,1) $ the sequence
$ \( (\One_{A_n}(x), \One_{B_n}(y) \)_{n=1}^\infty $ of elements of $ \{0,1\}
\times \{0,1\} $ is \normal{P(x,y)} in the following sense.

Given a probability measure $ p = (p_{00},p_{01},p_{10},p_{11}) $ on the
four-point set $ \{0,1\} \times \{0,1\} $ and a sequence of pairs $ \(
(a_n,b_n) \)_n $ where $ a_n,b_n \in \{0,1\} $, we say that the sequence is
\normal{p}, if it has appropriate frequencies of finite blocks, namely,
\[
\frac1n \#\{k\le n: (a_k,b_k) = (s_0,t_0), \dots, (a_{k+i},b_{k+i}) =
(s_i,t_i)\} \to p_{s_0,t_0} \dots p_{s_i,t_i}
\]
as $ n \to \infty $, for all $ i = 0,1,2,\dots $ and all $ s_0,\dots,s_i,
t_0,\dots,t_i \in\{0,1\} $.

The case $ i=0 $ gives $ \frac1n \#\{k\le n: x\in A_k, y\in B_k\} \to p_{11} $
and $ \frac1n \#\{k\le n: x\in A_k, y\notin B_k\} \to p_{10} $, thus,
\[
\frac1n \#\{k\le n: x\in A_k\} \to P_{10}(x,y) + P_{11}(x,y) \, ;
\]
note that $ P_{10}(x,y) + P_{11}(x,y) $ appears to be a function of $ x $
only. We see that $ f = P_{10} + P_{11} $, $ g = P_{01} + P_{11} $, $ h =
P_{11} $. It remains to prove that $ h = fg $.

We apply the bounded convergence theorem to the relation (integrated in $ x,y
$)
\begin{multline*}
\frac1n \#\{k\le n: (\One_{A_k}(x),\One_{B_k}(y)) = (s_0,t_0), \dots,
 ( \One_{A_{k+i}}(x),\One_{B_{k+i}}(y) ) = (s_i,t_i) \} \\
\to P_{s_0,t_0}(x,y) \dots P_{s_i,t_i}(x,y) \, ,
\end{multline*}
getting
\begin{multline*}
\frac1n \sum_{k=1}^n \mes_2 \{ (x,y) : \One_{A_k}(x) = s_0, \One_{B_k}(y) =
 t_0, \dots, \One_{A_{k+i}}(x) = s_i, \One_{B_{k+i}}(y) = t_i \} \\
\to \iint P_{s_0,t_0}(x,y) \dots P_{s_i,t_i}(x,y) \, \D x \D y \, .
\end{multline*}
(Here $ \mes_2 $ stands for the two-dimensional Lebesgue measure.) Summation
over $ t_0,\dots,t_i $ gives
\[
\frac1n \sum_{k=1}^n \mes \{ x : \One_{A_k}(x) = s_0, \dots, \One_{A_{k+i}}(x)
= s_i \} \to \int P_{s_0,*}(x) \dots P_{s_i,*}(x) \, \D x
\]
where $ P_{s,*} = P_{s,0} + P_{s,1} $; similarly,
\[
\frac1n \sum_{k=1}^n \mes \{ y : \One_{B_k}(y) = t_0, \dots, \One_{B_{k+i}}(y)
= t_i \} \to \int P_{*,t_0}(y) \dots P_{*,t_i}(y) \, \D y \, .
\]
Applying Ramsey's theorem (before Aldous' theorem) we ensure existence of
limits
\begin{gather*}
\lim_{k\to\infty} \mes \{ x : \One_{A_k}(x) = s_0, \dots, \One_{A_{k+i}}(x)
 = s_i \} \, , \\
\lim_{k\to\infty} \mes \{ y : \One_{B_k}(y) = t_0, \dots, \One_{B_{k+i}}(y) =
 t_i \} \, .
\end{gather*}
Taking into account that
\begin{multline*}
\mes_2 \{ (x,y) : \One_{A_k}(x) = s_0, \One_{B_k}(y) =
 t_0, \dots, \One_{A_{k+i}}(x) = s_i, \One_{B_{k+i}}(y) = t_i \} = \\
\mes \{ x : \One_{A_k}(x) = s_0, \dots, \One_{A_{k+i}}(x) = s_i \}
 \mes \{ y : \One_{B_k}(y) = t_0, \dots, \One_{B_{k+i}}(y) = t_i \}
\end{multline*}
we get
\begin{multline*}
\iint P_{s_0,t_0}(x,y) \dots P_{s_i,t_i}(x,y) \, \D x \D y = \\
= \bigg( \int P_{s_0,*}(x) \dots P_{s_i,*}(x) \, \D x \bigg)
 \bigg( \int P_{*,t_0}(y) \dots P_{*,t_i}(y) \, \D y \bigg)
\end{multline*}
for all $ i $ and $ s_0,t_0,\dots,s_i,t_i $. It means that
\begin{multline*}
\iint P_{00}^\al(x,y) P_{01}^\be(x,y) P_{10}^\ga(x,y) P_{11}^\de(x,y) \, \D x
 \D y = \\
= \bigg( \int P_{0*}^{\al+\be}(x) P_{1*}^{\ga+\de}(x) \, \D x \bigg)
  \bigg( \int P_{*0}^{\al+\ga}(y) P_{*1}^{\be+\de}(y) \, \D y \bigg)
\end{multline*}
for all $ \al,\be,\ga,\de \in \{0,1,2,\dots\} $. We see that all moments of
(the joint distribution of) $ P_{00}, P_{01}, P_{10}, P_{11} $ are equal to
moments of $ P_{0*}P_{*0}, P_{0*}P_{*1}, P_{1*}P_{*0}, P_{1*}P_{*1} $. Thus,
the distribution is the same. Therefore $ P_{s,t}(x,y) = P_{s,*}(x) P_{*,t}(y)
$ almost everywhere for $ s,t \in \{0,1\} $. In particular, $ P_{11} =
(P_{10}+P_{11}) (P_{01}+P_{11}) $, in other words, $ h = fg $.
\end{proof}

Recall that $ \limsup_n C_n = \cap_n (C_n \cup C_{n+1} \cup \dots ) $ for
arbitrary sets $ C_1, C_2, \dots $

\begin{proposition}\label{2.4}
For every measurable sets $ A_1, B_1, A_2, B_2, \dots \subset (0,1) $ there
exist measurable sets $ A,B \subset (0,1) $ such that
\begin{gather*}
\mes A \ge \liminf_n \mes A_n \, , \quad
 \mes B \ge \liminf_n \mes B_n \, , \\
A \times B \subset \limsup_n ( A_n \times B_n ) \, .
\end{gather*}
\end{proposition}

\begin{proof}
Taking into account the general relation $ \limsup_k C_{n_k} \subset \limsup_n
C_n $ we may replace the given sequence of pairs $ (A_n,B_n) $ with any
subsequence $ (A_{n_k},B_{n_k}) $ such that the limits
\[
a = \lim_k \mes A_{n_k} \, , \quad b = \lim_k \mes B_{n_k}
\]
exist. We forget the original sequence and rename the subsequence into $
(A_n,B_n) $. By Lemma \ref{2.3} we may assume that the limits
\begin{gather*}
f(x) = \lim_{n\to\infty} \frac{ \#\{k\le n: x\in A_k\} }{ n } \, , \quad
 g(y) = \lim_{n\to\infty} \frac{ \#\{k\le n: y\in B_k\} }{ n } \, , \\
h(x,y) = \lim_{n\to\infty} \frac{ \#\{k\le n: x\in A_k, y\in B_k\} }{ n }
\end{gather*}
exist almost everywhere, and $ h(x,y) = f(x) g(y) $ almost everywhere. Note
that
\[
\int f(x) \, \D x = \lim_n \frac1n (\mes A_1 + \dots + \mes A_n) = a \, ,
\]
and similarly $ \int g(y) \, \D y = b $.

We have
\[
\limsup_n (A_n \times B_n) \supset \{ (x,y) : h(x,y) > 0 \}
\]
just because a set of nonzero frequency cannot be finite. Taking
\[
A = \{ x : f(x) > 0 \} \, , \quad B = \{ y : g(y) > 0 \} \, ,
\]
we get
\[
 \limsup_n ( A_n \times B_n ) \supset A \times B
\]
since  $ h(x,y) = f(x) g(y) > 0 $ for $ (x,y) \in A \times B $. It remains to
check that $ \mes A \ge a $ and $ \mes B \ge b $, which is easy:
\[
a = \int f(x) \, \D x = \int_A f(x) \, \D x \le \mes A
\]
and similarly for $ B $.
\end{proof}

\begin{proof}[Proof of Prop.~\textup{\ref{2.2}}]
We choose $ U_n, V_n $ such that $ \( U_n \times (0,1) \) \cup \( (0,1) \times
V_n \) \supset W_n $ and $ \mes(U_n) + \mes(V_n) \to \lim_n \be(W_n) $. Taking
a subsequence (which does not change $ \lim_n \be(W_n) $) we ensure existence
of the limits $ \lim_n \mes(U_n) $, $ \lim_n \mes(V_n) $. The complements $
A_n = (0,1) \setminus U_n $, $ B_n = (0,1) \setminus V_n $ satisfy
\begin{gather*}
(A_n \times B_n) \cap W_n = \emptyset \, , \\
\mes(A_n) \to a \, , \quad \mes(B_n) \to b \quad \text{as } n \to \infty \, ,
 \\
(1-a) + (1-b) = \lim_n \be(W_n) \, .
\end{gather*}
Prop.~\ref{2.4} gives us $ A,B $ such that $ \mes A \ge a $, $ \mes B \ge b $
and $ A \times B \subset \limsup_n (A_n \times B_n) $. Taking the complements
$ U = (0,1) \setminus A $, $ V = (0,1) \setminus B $ we get $ \mes U + \mes V
\le \lim_n \be(W_n) $ and
\begin{multline*}
\( U \times (0,1) \) \cup \( (0,1) \times V \) \supset \liminf_n \( U_n \times
 (0,1) \) \cup \( (0,1) \times V_n \) \supset \\
\supset \liminf_n W_n = W \, ,
\end{multline*}
therefore $ \be(W) \le \lim_n \be(W_n) $. The converse inequality is trivial.
\end{proof}

The proof gives us the following by-product.

\begin{corollary}
The infimum in the definition of $ \be(W) $ is always reached.
\end{corollary}

(The special case $ \be(W)=0 $ is much simpler, see \cite[5.7]{I}.)

I do not know whether the supremum in the definition of $ \al(W) $ is reached
or not in general. However it is reached in the special case stated below
(only this case is used in the next section).

\begin{proposition}\label{2.6}
Let $ W \subset (0,1) \times (0,1) $ be a set of class $ \A_{\de\si} $ such
that $ W \cap (A \times B) \ne \emptyset $ for all measurable sets $ A,B
\subset (0,1) $ of positive measure. Then there exists $ m \in M $ such that $
m(W) = 1 $.
\end{proposition}

\begin{proof}
If $ W \subset \( U \times (0,1) \) \cup \( (0,1) \times V \) $ then $ \mes U
= 1 $ or $ \mes V = 1 $; therefore $ \be(W) = 1 $. By Prop.~\ref{2.2}, $
\al(W) = \be(W) = 1 $. It remains to prove that the supremum in the definition
of $ \al(W) $ is reached.

As was noted after Prop.~\ref{2.2}, its true generality is not restricted to $
(0,1) $ with Lebesgue measure. Especially, if $ \mu,\nu $ are absolutely
continuous probability measures on $ (0,1) $ then clearly $ W \cap (A \times
B) \ne \emptyset $ whenever $ \mu(A) > 0 $, $ \nu(B) > 0 $; a generalization
of Prop.~\ref{2.2} gives $ \al_{\mu,\nu} (W) = 1 $ where $ \al_{\mu,\nu} (W) =
\sup \{ m(W) : m \in M_{\mu,\nu} \} $ and $ M_{\mu,\nu} $ is the set of all
joinings between $ \mu $ and $ \nu $.

The set $ M_W = \{ m \in M : m \( (0,1) \setminus W \) = 0 \} $ contains the
limit of each \emph{increasing} sequence of elements of $ M_W $. It follows
easily that $ M_W $ contains a maximal element $ m $ (at least one). It
remains to prove that $ m $ is a probability measure.
Assume the contrary: $ m \( (0,1) \times (0,1) \) < 1 $. Consider the
marginals $ m_1, m_2 $ of $ m $ and the measures $ \mu = \mes - m_1 $, $ \nu =
\mes - m_2 $ (where `$ \mes $' is the Lebesgue measure on $ (0,1) $); $ \mu $
and $ \nu $ are positive measures on $ (0,1) $, absolutely continuous, and $
\mu\((0,1)\) = \nu\((0,1)\) > 0 $. According to the previous paragraph
(applied to normalized $ \mu,\nu $) there exists a positive measure $ n $ on $
W $ whose marginals $ n_1, n_2 $ satisfy $ n_1 \le \mu $, $ n_2 \le \nu $ and
such that $ n(W) $ is close to $ \mu\((0,1)\) $; however, we only need to know
that $ n(W) > 0 $. The measure $ m+n $ belongs to $ M $, which contradicts to
the maximality of $ m $.
\end{proof}

\section[]{\raggedright Selectors and independence}
\label{sect3}\begin{definition}\label{3.1}
Let $ X : \Om \to \DCS(0,1) $ be a random dense countable subset of $ (0,1)
$. A \emph{selector} (of $ X $) is a random variable $ Z : \Om \to (0,1) $
such that
\[
Z(\om) \in X(\om) \quad \text{for almost all } \om \, .
\]
\end{definition}

In terms of a chosen measurable enumeration $ X = \{ Y_1, Y_2, \dots \} $,
the general form of a selector is
\[
Z(\om) = Y_{n(\om)} (\om) \quad \text{for almost all } \om \, ,
\]
where $ n $ runs over all measurable maps $ \Om \to \{ 1,2,\dots \} $.

Sometimes dependence between two random variables reduces to a joint density
(w.r.t.\ their marginal distributions). Here are two formulation in general
terms.

\begin{lemma}\label{3.2}
Let $ (\Om,\F,P) $ be a probability space and $ C \subset \Om $ a measurable
set. The following two conditions on a pair of sub-\sif s $ \F_1, \F_2 \subset
\F $ are equivalent:

(a) there exists a measurable function $ f : \Om \times \Om \to [0,\infty) $
such that
\[
P ( A \cap B \cap C ) = \int_{A\times B} f(\om_1,\om_2) \, P(\D\om_1)
P(\D\om_2) \quad \text{for all } A \in \F_1, B \in \F_2 \, ;
\]

(b) there exists a measurable function $ g : C \times C \to [0,\infty) $
such that
\[
P ( A \cap B \cap C ) = \int_{(A\cap C)\times(B\cap C)} g(\om_1,\om_2) \,
P(\D\om_1) P(\D\om_2) \quad \text{for all } A \in \F_1, B \in \F_2 \, .
\]
\end{lemma}

(Note that $ f,g $ may vanish somewhere.)

\begin{proof}
(b) \imp (a): just take $ f(\om_1,\om_2) = g(\om_1,\om_2) $ for $ \om_1,\om_2
  \in C $ and $ 0 $ otherwise.

(a) \imp (b):
we consider conditional probabilities $ h_1 = \cP{C}{\F_1} $, $ h_2 =
\cP{C}{\F_2} $, note that $ h_1(\om) > 0 $, $ h_2(\om) > 0 $ for almost all $
\om \in C $ and define
\[
g(\om_1,\om_2) = \frac{ f(\om_1,\om_2) }{ h_1(\om_1) h_2(\om_2) } \quad
\text{for } \om_1, \om_2 \in C \, .
\]
Then
\begin{multline*}
\int_{(A\cap C)\times(B\cap C)} g(\om_1,\om_2) \, P(\D\om_1) P(\D\om_2) = \\
= \int_{A\times B} \frac{ f(\om_1,\om_2) \One_C(\om_1) \One_C(\om_2) }{
 h_1(\om_1) h_2(\om_2) } \, P(\D\om_1) P(\D\om_2) \, .
\end{multline*}
(The integrand is treated as $ 0 $ outside $ C \times C $.)
Assuming that $ f $ is \measurable{(\F_1\otimes\F_2)} (otherwise $ f $ may be
replaced with its conditional expectation) we see that the conditional
expectation of the integrand, given $ \F_1\otimes\F_2 $, is equal to $
f(\om_1,\om_2) $. Thus, the integral is
\[
\dots = \int_{A\times B} f(\om_1,\om_2) \, P(\D\om_1) P(\D\om_2) = P ( A \cap
B \cap C ) \, .
\]
\end{proof}

\begin{definition}\label{3.3}
Let $ (\Om,\F,P) $ be a probability space and $ C \subset \Om $ a measurable
set. Two sub-\sif s $ \F_1, \F_2 \subset \F $ are a \emph{nonsingular pair}
within $ C $, if they satisfy the equivalent conditions of Lemma \ref{3.2}.
\end{definition}

\begin{lemma}\label{3.4}
(a)
Let $ C_1 \subset C_2 $. If $ \F_1, \F_2 $ are a nonsingular pair within $ C_2
$ then they are a nonsingular pair within $ C_1 $.

(b)
Let $ C_1,C_2,\dots $ be pairwise disjoint and $ C = C_1 \cup C_2 \cup \dots $
If $ \F_1, \F_2 $ are a nonsingular pair within $ C_k $ for each $ k $ then
they are a nonsingular pair within $ C $.

(c)
Let $ \Ec_1 \subset \F $ be another sub-\sif\ such that $ \Ec_1 \subset \F_1 $
within $ C $ in the sense that
\[
\forall E \in \Ec_1 \;\; \exists A \in \F_1 \;\; ( A \cap C = E \cap C ) \, .
\]
If $ \F_1, \F_2 $ are a nonsingular pair within $ C $ then $ \Ec_1, \F_2 $
are a nonsingular pair within $ C $.
\end{lemma}

\begin{proof}
(a)
We define two measures $ \mu_1, \mu_2 $ on $ (\Om,\F_1) \times (\Om,\F_2) $ by
$ \mu_k (Z) = P ( C_k \cap \{ \om : (\om,\om) \in Z \} ) $ for $ k=1,2
$. Clearly, $ \mu_k (A\times B) = P ( A \cap B \cap C_k ) $. Condition
\ref{3.2}(a) for $ C_k $ means absolute continuity of $ \mu_k $ (w.r.t.\ $ P
\times P $). However, $ \mu_1 \le \mu_2 $.

(b)
Using the first definition, \ref{3.2}(a), we just take $ f = f_1 + f_2 + \dots
$

(c)
Immediate, provided that the second definition us used, \ref{3.2}(b).
\end{proof}

\begin{proposition}\label{3.5}
Let $ X $ be a random dense countable subset of $ (0,1) $ such that two random
sets $ X \cap (0,\frac12) $, $ X \cap (\frac12,1) $ are independent. Let $
Z_1,Z_2,\dots $ be a measurable enumeration of $ X \cap (\frac12,1) $
independent of \emph{some} measurable enumeration of $ X \cap (0,\frac12)
$. Let $ Y_1,\dots,Y_n $ be selectors of $ X $. Then the sub-\sif\ $ \F_1 $
generated by $ Y_1,\dots,Y_n $ and the sub-\sif\ $ \F_2 $ generated by $
Z_1,Z_2,\dots $ are a nonsingular pair within the event
\[
C = \{ Y_1,\dots,Y_n < \tfrac12 \} \, .
\]
\end{proposition}

The proof is given below after some discussion.
The condition imposed on $ X $ is the relevant part of the independence
condition (recall Def.~\ref{1.2}). The threshold is chosen at $ \frac12 $, but
any other number of $ (0,1) $ could be used equally well.

It may happen that $ \Pr{C} = 1 $; even in this case $ \F_1, \F_2 $ need not
be independent (see Counterexample \ref{3.6}). Of course, there exist
an enumeration $ X \cap (0,\frac12) = \{ U_1, U_2, \dots \} $ independent of
$ (Z_k)_k $ and, for instance, $ Y_1(\om) = U_{n(\om)}(\om) $. However, $
n(\cdot) $ need not be independent of $ Z_1, Z_2, \dots $

Note the finite sequence $ (Y_k)_k $ but the infinite sequence $ (Z_k)_k
$. The proposition may fail for an infinite sequence $ (Y_k)_k $ (see
Counterexample \ref{3.6}). Note also that $ Y_k $ need not be pairwise
different.

\begin{counterexample}\label{3.6}
We start with independent random variables $ U_k, V_k, Z_k $ ($ k = 1,2,\dots
$) such that each $ U_k $ is distributed uniformly on $ (0,\frac14) $, each $
V_k $ --- on $ (\frac14, \frac12) $, and each $ Z_k $ --- on $ (\frac12,1)
$. We construct events $ A_1, A_2, \dots $ that generate the sub-\sif\ $ \F_2
$ generated by $ Z_1,Z_2,\dots $ Now we define random variables $ Y_1, Y_2,
\dots $ as follows:
\[
Y_k = \begin{cases}
 U_k &\text{on $ A_k $},\\
 V_k &\text{outside $ A_k $}.
\end{cases}
\]
Clearly, each $ Y_k $ is a selector of $ X = \{U_1,U_2,\dots\} \cup
\{V_1,V_2,\dots\} \cup \{Z_1,Z_2,\dots\} $, and $ A_k = \{ Y_k < \frac14 \}
$. The sub-\sif\ $ \F_1 $ generated by $ Y_1,Y_2,\dots $ contains $
A_1,A_2,\dots $, therefore
\[
\F_1 \supset \F_2 \, .
\]
Clearly, $ \F_1 $ and $ \F_2 $ are not a nonsingular pair.
\end{counterexample}

\begin{proof}[Proof of Prop.~\textup{\ref{3.5}}]
I assume that $ n=1 $ (thus, $ Y=Y_1 $), leaving to the reader the
straightforward generalization. We choose a measurable enumeration $ U_1, U_2,
\dots $ of $ X \cap (0,\frac12) $ independent of $ Z_1,Z_2,\dots $ and
partition $ C $ into $ C_k = \{ Y = U_k \} $. By \ref{3.4}(b) it is sufficient
to prove that $ \F_1, \F_2 $ are a nonsingular pair within each $ C_k $. By
\ref{3.4}(c) we replace $ \F_1 $ with the \sif\ $ \si(U_k) $ generated by $
U_k $. By \ref{3.4}(a) we replace $ C_k $ with the whole $ \Om $. The \sif s $
\si(U_k), \F_2 $ are a nonsingular pair within $ \Om $, since they are
independent.
\end{proof}

\section[]{\raggedright Existence of selectors}
\label{sect4}In order to allow for some additional randomization, in this section we often
construct a selector not on the original probability space $ \Om $ but on some
extended space $ \ti\Om $. In fact, the product space $ \Om \times \Om' $
(with the product measure), where $ \Om' $ is a nonatomic probability space,
may serve as $ \ti\Om $. Naturally, the given random dense countable set $ X $
is transferred to $ \ti\Om $ by $ X(\om,\om') = X(\om) $.

\begin{lemma}\label{4.1}
The following two conditions on a random dense countable set $ X : \Om \to
\DCS(0,1) $ are equivalent:

(a)
there exists a selector $ Z : \ti\Om \to (0,1) $ of $ X $ (on some extension $
\ti\Om $ of $ \Om $), distributed uniformly on $ (0,1) $;

(b)
there exists a probability measure on the set
\begin{equation}\label{4.2}
W = \{ (\om,t) : t \in X(\om) \} \subset \Om \times (0,1)
\end{equation}
whose two marginals are $ P $ and the uniform distribution on $ (0,1) $.
\end{lemma}

\begin{proof}
(a) \imp (b):
The joint distribution of $ \om $ (treated as a function of $ \ti\om \in
\ti\Om $) and $ Z $ is the needed measure on $ W $.

(b) \imp (a):
We take $ \ti\Om = \Om \times \Om' $, disintegrate the given measure $ m $ on
$ W $ into conditional measures $ m_\om $ on $ (0,1) $, represent $ m_\om $ as
the distribution of some $ Z_\om : \Om' \to (0,1) $ and combine these $ Z_\om
$ by $ Z(\om,\om') = Z_\om (\om') $. Measurability of $ Z $ can be achieved by
choosing $ \Om' $ to be $ (0,1) $ (with Lebesgue measure) and each $ Z_\om $
to be an \emph{increasing} function $ (0,1) \to (0,1) $.
\end{proof}

The set $ W $ defined by \eqref{4.2} should be treated modulo sets of the form
$ A \times (0,1) $, $ P(A) = 0 $ (and their subsets). Then $ W $ appears to
belong to the class $ \A_{\de\si} $ (introduced before Lemma \ref{2.3}), as
stated below.

\begin{lemma}\label{4.3}
Let $ X : \Om \to \DCS(0,1) $ be a random dense countable set. Then there
exists a subset $ \Om_1 \subset \Om $ of probability $ 1 $ such that the set
\[
W_1 = \{ (\om,t) : \om \in \Om_1, t \in X(\om) \}
\]
belongs to the class $ \A_{\de\si} $.
\end{lemma}

\begin{proof}
We take a measurable enumeration $ X = \{ Y_1, Y_2, \dots \} $ and note that $
W = W_1 \cup W_2 \cup \dots $ where $ W_k = \{ (\om,Y_k(\om)) : \om\in\Om \} $
is the graph of $ Y_k $. We choose a topology on $ \Om $ turning $ \Om $ into
a compact metrizable space (such that the given probability measure on $ \Om $
is a Borel measure, up to negligible sets). By Lusin's theorem there exist
compact sets $ C_n \subset \Om $ such that the restrictions $ Y_k |_{C_n} $
are continuous and the set $ \Om_1 = C_1 \cup C_2 \cup \dots $ is of
probability $ 1 $. The \emph{compact} sets $ W_{k,n} = \{ (\om,Y_k(\om) :
\om\in C_n \} $ belong to the class $ \A_\de $, therefore their union $ W_1 $
belongs to $ \A_{\de\si} $.
\end{proof}

\begin{proposition}\label{4.4}
Let $ X : \Om \to \DCS(0,1) $ satisfy \eqref{1.7}. Then $ X $ has a selector $
Z : \ti\Om \to (0,1) $ (on some extension $ \ti\Om $ of $ \Om $), distributed
uniformly on $ (0,1) $.
\end{proposition}

\begin{proof}
Lemma \ref{4.3} gives us a set $ W_1 \subset W $ of class $ \A_{\de\si} $;
\eqref{1.7} shows that $ W_1 $ intersects every $ A \times B $ where $ A
\subset \Om $, $ P(A) > 0 $ and $ B \subset (0,1) $, $ \mes B > 0
$. Proposition \ref{2.6} (or rather, its evident generalization) gives us a
measure on $ W_1 $ whose marginals are $ P $ and $ \mes $. Lemma \ref{4.1}
((b) \imp (a)) completes the proof.
\end{proof}

Only the second part of \eqref{1.7} was used.

\begin{lemma}\label{4.5}
Let $ X $ be a random dense countable subset of $ (0,1) $. If for every $ \eps
> 0 $ the random dense countable subset $ X \cap (\eps,1) $ of $ (\eps,1) $
satisfies \eqref{1.7}, then $ X $ satisfies \eqref{1.7}.
\end{lemma}

The (straightforward) proof is left to the reader. More generally, one may
check $ X \cap (0,\frac12-\eps) $ and $ X \cap (\frac12+\eps,1) $, etc.

\begin{proposition}\label{4.6}
Let $ X $ be a random dense countable subset of $ (0,1) $, satisfying
\eqref{1.7} and the independence condition. Let $ Y_1, \dots, Y_n $ be
selectors of $ X $. Then there exists a selector $ Z $ (on some extension of
the given probability space) independent of $ Y_1, \dots, Y_n $ and
distributed uniformly on $ (0,1) $.
\end{proposition}

\begin{proof}
I assume that $ n=1 $ (thus, $ Y = Y_1 $), leaving to the reader the
straightforward generalization. We disintegrate the given measure $ P $ on the
given probability space $ (\Om,\F,P) $ into conditional (given $ Y=y $)
measures $ P_y $ on $ \Om $ for $ y \in (0,1) $; $ P = \int P_y \, \mu(\D y)
$, where $ \mu $ is the distribution of $ Y $.

Being considered w.r.t.\ $ P_y $, the \measurable{\F} map $ X : \Om \to
\DCS(0,1) $ is another random dense countable set; denote it $ X_y $. We claim
that
\begin{equation}\label{4.7}
X_y \text{ satisfies \eqref{1.7} for \almost{\mu} every } y \in (0,1) \, .
\end{equation}
By (generalized) Lemma \ref{4.5} it is enough to check \eqref{1.7} for $ X_y
\cap (0,a) $ for \almost{\mu} all $ y \in (a,1) $ as well as $ X_y \cap (a,1)
$, $ y \in (0,a) $; here $ a $ runs over $ (0,1) $ (or only its rational
points). Proposition \ref{3.5} (generalized a bit) shows that, roughly
speaking, the distribution of $ X_y \cap (a,1) $ is absolutely continuous
w.r.t.\ the distribution of $ X \cap (a,1) $, as far as $ y \in (0,a) $. Thus,
\eqref{1.7} for $ X \cap (a,1) $ implies \eqref{1.7} for $ X_y \cap (a,1) $,
and \eqref{4.7} is verified.

Proposition \ref{4.4} gives us selectors $ Z_y : \ti\Om \to (0,1) $ of $ X $
such that $ Z_y $ is uniformly distributed on $ (\Om,P_y) $. We want to glue
them together,
\[
Z(\om) = Z_{Y(\om)} (\om) \, ;
\]
to this end, however, we need \measurability{\mu} of $ Z_y $ in $ y $, which
is the goal below.

Similarly to the proof of Lemma \ref{4.3} we take a measurable enumeration $ X
= \{ U_1, U_2, \dots \} $, turn $ \Om $ into a compact metrizable space and
construct compact sets $ C_n \subset \Om $ such that the restrictions $ U_k
|_{C_n} $ are continuous and the set $ \Om_1 = C_1 \cup C_2 \cup \dots $
satisfies $ P(\Om_1) = 1 $. Then $ P_y (\Om_1) = 1 $ for \almost{\mu} all $ y
$ (since $ P = \int P_y \, \mu(\D y) $). Once again, the \set{\A_{\de\si}} $
W_1 = \{ (\om,t) : \om \in \Om_1, t \in X(\om) \} $ is the union of compact
sets $ W_{k,n} = \{ (\om,U_k(\om)) : \om \in C_n \} $. The set $ M $ of all
probability measures on $ W_1 $ is a standard Borel space. The same holds for
measures on $ \Om $ and on $ (0,1) $. Denote the two marginals of a measure $
m \in M $ by $ \phi(m) $ and $ \psi(m) $; $ \phi,\psi $ are Borel measurable
maps. We know that the set
\[
M_y = \{ m \in M : \phi(m) = P_y, \, \psi(m) = \mes \}
\]
is nonempty for \almost{\mu} all $ y $. (Indeed, the selectors $ Z_y $ are
constructed in the proof of Prop.~\ref{4.4} via measures that belong to $ M_y
$.) Taking into account that $ P_y $ is (or rather, may be chosen to be) a
Borel measurable function of $ y $, we apply well-known uniformization
theorems and get a \measurable{\mu} map $ y \mapsto m_y $ such that $ m_y \in
M_y $ for \almost{\mu} all $ y $. Now \measurability{\mu} of the map $ y
\mapsto Z_y $ can be achieved (recall the end of the proof of Lemma
\ref{4.1}).
\end{proof}

\section[]{\raggedright Proving the main results}
\label{sect5}\begin{proof}[Proof of Theorem \textup{\ref{1.8}}]
Let $ X $ be a random dense countable subset of $ (0,1) $ satisfying
\eqref{1.7} and the independence condition. In order to prove that $ X $ is
distributed like the unordered infinite sample it is sufficient to find a
measurable enumeration $ X = \{ Y_1, Y_2, \dots \} $ that satisfies the
conditions of the main lemma (Lemma 2.1) of \cite{I}. Namely, the conditional
distribution of $ Y_n $ given $ Y_1,\dots,Y_{n-1} $ must have a density $
(t,\om) \mapsto f_n(t,\om) $ and the series $ \sum_{n=1}^\infty f_n(t,\om) $
must diverge for almost all pairs $ (t,\om) $.

Existence of conditional densities is ensured by \eqref{4.7} (generalized for
$ n \ge 1 $) for every measurable enumeration.

Starting with an arbitrary measurable enumeration $ (Z_k)_k $ we construct the
needed enumeration $ (Y_k)_k $ as follows. First, $ Y_1 = Z_1 $. Second,
Prop.~\ref{4.6} gives us $ Y_2 $ independent of $ Y_1 $ and distributed
uniformly on $ (0,1) $. (The probability space is extended as needed.) Third,
$ Y_3 = Z_2 $ unless $ Y_2 = Z_2 $, in which case $ Y_3 = Z_3
$. Prop.~\ref{4.6} gives us $ Y_4 $ independent of $ Y_1, Y_2, Y_3 $ and
uniform. And so on; $ Y_{2n} $ is independent of $ Y_1,\dots,Y_{2n-1} $ and
uniform, while $ Y_{2n+1} $ is the first $ Z_k $ different from $
Y_1,\dots,Y_{2n} $. Clearly, $ (Y_k)_k $ is a measurable enumeration of $ X $,
and $ \sum f_n = \infty $ since $ f_{2n}(t,\om) = 1 $.
\end{proof}

It remains to prove Proposition \ref{1.10}.

Recall the cyclic shift $ T_s $ (introduced after \ref{1.8}).

\begin{lemma}\label{5.2}
Let $ A \subset (0,1) $ be a measurable set, $ \mes A > 0 $, and $ L \subset
(0,1) $ a dense countable set (not random). Then the set $ A \cap T_s(L) $ is
infinite for almost all $ s \in (0,1) $.
\end{lemma}

\begin{proof}
First, the set $ B = \{ s : A \cap T_s(L) \ne \emptyset \} = \cup_{l\in L}
T_l^{-1} (A) $ is of full measure, since $ \frac1\eps \mes \( A \cap
(a,a+\eps) \) \le \frac1\eps \mes \( B \cap (b,b+\eps) \) $ for all $ \eps > 0
$ and $ a,b \in (0,1-\eps) $.

Second, $ L \supset L_1 \cup L_2 \cup \dots $ for some pairwise disjoint dense
countable sets $ L_1, L_2, \dots $ Almost every $ s $ satisfies $ A \cap
T_s(L_n) \ne \emptyset $ for all $ n $.
\end{proof}

\begin{proof}[Proof of Proposition \textup{\ref{1.10}}]
Condition \eqref{1.7} is satisfied by the given $ X : \Om \to \DCS(0,1) $ if
(and only if) it is satisfied by $ Y : \Om \times (0,1) \to \DCS(0,1) $
used in \ref{1.9}, $ Y(\om,s) = T_s(X(\om)) $ (just because $ X $ and $ Y $
are identically distributed). If $ \mes A = 0 $ then $ A \cap Y(\om,s) =
\emptyset $ for all $ s $ except for the negligible set $ \cup_{l\in X(\om)}
T_l^{-1}(A) $. If $ \mes A > 0 $ then $ A \cap Y(\om,s) $ is infinite for
almost all $ s $ by Lemma \ref{5.2}, and it holds for almost all $ \om $.
\end{proof}

\bigskip
\filbreak
{
\small
\begin{sc}
\parindent=0pt\baselineskip=12pt
\parbox{4in}{
Boris Tsirelson\\
School of Mathematics\\
Tel Aviv University\\
Tel Aviv 69978, Israel
\smallskip
\par\quad\href{mailto:tsirel@post.tau.ac.il}{\tt
 mailto:tsirel@post.tau.ac.il}
\par\quad\href{http://www.tau.ac.il/~tsirel/}{\tt
 http://www.tau.ac.il/\textasciitilde tsirel/}
}

\end{sc}
}
\filbreak

\end{document}